\documentclass{amsart}
\usepackage{amsmath}
\usepackage{graphicx}
\usepackage{latexsym}

\newtheorem{thm}{Theorem}[section]
\newtheorem{lem}[thm]{Lemma}
\newtheorem{cor}[thm]{Corollary}
\newtheorem{prop}[thm]{Proposition}

\theoremstyle{definition}
\newtheorem{defn}[thm]{Definition}

\newcommand{\R}{\mathbb{R}}

\begin{document}

\title[Cords and $1$-handles attached to surface-knots]
{Cords and $1$-handles attached to surface-knots}
\vspace{0.2in} 

\author[S. Kamada]{Seiichi Kamada}
\address{{Department of Mathematics, Osaka City University} \\ 
{Sumiyoshi, Osaka 558-8585, Japan}} 
\email{skamada{@}sci.osaka-cu.ac.jp}

\begin{abstract}
J. Boyle classified $1$-handles attached to surface-knots, that are closed and connected surfaces embedded in the Euclidean $4$-space, in the case that the surfaces are oriented and  $1$-handles are orientable with respect to the orientations of the surfaces.  In that case, 
the equivalence classes of $1$-handles correspond to the equivalence classes of cords attached to the surface-knot, and correspond to 
the double cosets of the peripheral subgroup of the knot group.  In this paper, we classify cords and cords with local orientations  attached to  (possibly non-orientable) surface-knots.  And we classify  $1$-handles attached to surface-knots in the case that the surface-knots are oriented and $1$-handles are non-orientable, and in the case that the surface-knots are non-orientable.   
\end{abstract}

\maketitle

\renewcommand{\thefootnote}{\fnsymbol{footnote}}
\setcounter{footnote}{-1}
\renewcommand{\thefootnote}{\arabic{footnote}}

\begin{center}
Dedicated to Professor Francisco Gonz{\' a}lez-Acu{\~n}a on his seventieth birthday
\end{center}

\section{Introduction}

By a {\it surface-knot} we mean a closed (possibly non-orientable) and connected surface embedded in $\R^4$.  
F. Hosakawa and A. Kawauchi \cite{HK} studied unknotted surface-knots in $\R^4$ and surgery along $1$-handles attached to surface-knots.   They proved that an oriented surface-knot $F$ in $\R^4$ satisfies that the knot group $\pi_1(\R^4 - F)$  is infinite cyclic if and only if an unknotted surface-knot can be obtained from $F$ by surgery along trivial $1$-handles.  
A similar result holds for a non-orientable surface-knot  in $\R^4$ (cf. \cite{Kam89}).   Surgery along a $1$-handle is a method of constructing a surface-knot  from another with lower genus.  The knot type of the surface-knot obtained from a surface-knot  in $\R^4$ by surgery along a $1$-handle depends on the equivalence class of the $1$-handle.  Classifying  $1$-handles attached to a surface-knot $F$ is important in order to consider the knot types obtained from $F$ by surgery along $1$-handles.  J. Boyle \cite{Boy} classified such $1$-handles in the case that $F$ is oriented, and $1$-handles are orientable with respect to the orientation of $F$.   This case, say (Case 1),  is sufficient when we work on oriented surface-knots in $R^4$.  When we work on non-orientable surface-knots, we should also consider  the following two cases: (Case 2) $F$ is oriented and $1$-handles are non-orientable with respect to the orientation of $F$, and (Case 3) $F$ is non-orientable.  In this paper we give a classification theorem to each of these two cases (Case 2) and (Case 3), which is analogous to Boyle's classification in (Case 1).  

In order to classify $1$-handles attached to a surface-knot, we first classify cords and cords with local orientations attached to a surface-knot.  Roughly speaking, the equivalence classes of $1$-handles attached to a surface-knot $F$ correspond to the equivalence classes of cords attached to $F$ in Cases 1 and 2, or correspond to the equivalence classes of cords with local orientations at the endpoints attached to $F$    in Case 3.  

We work in the PL category and all embedded surfaces in $4$-manifolds are assumed to be locally flat. The results in this paper are also valid in the smooth category.  

Throughout this paper, $B^n$ denotes the unit $n$-ball in $\R^n$ and $0 \in B^n$ is the center.  

An {\it ambient isotopy} of a space $X$ is 
an isotopy $(f_s \,| \, s \in [0,1])$ such that for each $s\in [0,1]$, $f_s: X \to X$ is a homeomorphism and 
$f_0$ is the identity map of $X$. Two subsets $A$ and $A'$ of $X$ are {\it ambient isotopic} if there is an ambient isotopy $(f_s \,| \, s \in [0,1])$ of $X$ with $f_1(A)=A'$.  Two maps $g: Y \to X$ and $g': Y \to X$ are {\it ambient isotopic} if there is an ambient isotopy $(f_s \,| \, s \in [0,1])$ of $X$ with $f_1 \circ g = g'$.

Some results of this paper are given partially in Section 5.2 of \cite{KamBook2}, written in Japanese.  This paper completes it. 

\section{Definitions on $1$-handles}

There are two notions of $1$-handles, one is defined by embeddings (cf. \cite{Boy}) and the other is defined by 3-cells in $\R^4$ (cf. \cite{HK}).   To distinguish these two, we call a $1$-handle as an embedding a {\it $1$-handle map} in this paper.  

Let $F$ be a surface-knot.  

\begin{defn}
A {\it $1$-handle map} attached to $F$ is an embedding $h: [0,1] \times B^2 \to \R^4$ with $F \cap h([0,1] \times B^2) = h(\{0,1\} \times B^2)$.  The restriction of $h$ to $[0,1] \times \{0\}$ $(=[0,1])$  is denoted by $h^c: [0,1] \to R^4$ and 
called the {\it core map}.  The image of $h^c$ is called the {\it core} of $h$.  
\end{defn}

For a $1$-handle map $h: [0,1] \times B^2 \to \R^4$, the {\it reverse} of $h$ is a $1$-handle map 
${\rm rev}(h) : [0,1] \times B^2 \to \R^4$ with ${\rm rev}(h)(t, x) = h(1-t, x)$.  

\begin{defn}
Let $h: [0,1] \times B^2 \to \R^4$ and $h': [0,1] \times B^2 \to \R^4$ be $1$-handle maps attached to $F$.  
\begin{itemize}
\item[(1)] $h$ and $h'$ are {\it equivalent} if they are ambient isotopic in $\R^4$ by an ambient isotopy of $\R^4$ keeping $F$ setwise fixed.  
\item[(2)] $h$ and $h'$ are {\it equivalent up to reversion} if $h$ is equivalent to $h'$ or ${\rm rev}(h')$. 
\end{itemize}
\end{defn}

For a $1$-handle map $h$ attached to $F$, we denote by 
${\rm h}^1(F; h)$ the surface-knot 
$$ (F - h( \{0,1\} \times {\rm int } B^2 ) \cup h ( [0,1] \times \partial B^2),$$
which we call the surface-knot obtained from $F$ by {\it surgery} along $h$.  The surgery is also called a $1$-handle surgery or a hyperboloidal transformation (\cite{HK}).  The symbol ${\rm h}^1$ stands for a $1$-handle surgery.  In \cite{Boy} it is denoted by $F + h$.  

If $h$ and $h'$ are equivalent or equivalent up to reversion attached to a surface-knot $F$, then 
${\rm h}^1(F; h)$ and ${\rm h}^1(F; h')$ are ambient isotopic in $\R^4$.  

\begin{defn}
Assume that $F$ is an orientable surface-knot.  A $1$-handle map $h$ attached to $F$ is {\it orientable} (or {\it non-orientable}, resp.)  if ${\rm h}^1(F; h)$ is orientable (or non-orientable, resp.).     
\end{defn}

When $F$ is oriented and $h$ is orientable, the surface-knot ${\rm h}^1(F; h)$ is assumed to have an orientation that coincides, over $ F - h( \{0,1\} \times {\rm int } B^2 )$, with the orientation of $F$.  

Now we recall the notion of a $1$-handle as a $3$-cell in $\R^4$ from \cite{HK}.  

\begin{defn}
A {\it $1$-handle} attached to $F$ is a 
3-cell $B$ in $\R^4$ such that $B \cap F = \partial B \cap F$ and this intersection is the union of disjoint two $2$-cells.  A properly embedded arc $C$ in $B$ is called a {\it core} of $B$ if it is a strong deformation retract of $B$ and it connects an interior point of one 2-cell of $B \cap F$ with another interior point of the other 2-cell.  
\end{defn}

For a $1$-handle map $h$ attached to $F$, the image of $h$ is a $1$-handle attached to $F$, say $B$, and the core of $h$ is a core of $B$.  Conversely, for a $1$-handle $B$ and a core $C$ of $B$, there is a $1$-handle map $h$ whose image is $B$ and its core is $C$. 

\begin{defn}
Two $1$-handles $B$ and $B'$ attached to $F$ are {\it equivalent} if they are ambient isotopic in $\R^4$ by an ambient isotopy of $\R^4$ keeping $F$ setwise fixed.  
\end{defn}

\begin{lem}\label{lem:twonotionsA}
For $1$-handles $B$ and $B'$ attached to $F$, let $h$ and $h'$ be $1$-handle maps attached to $F$ whose images are $B$ and $B'$, respectively.  $B$ and $B'$ are equivalent if and only if $h$ and $h'$ are equivalent up to reversion.  
\end{lem}

This lemma follows from Lemma~\ref{lem:twonotionsB} stated below.

\begin{defn}
A {\it $1$-handle with an oriented core} attached to $F$ is a pair $(B,C)$ of a $1$-handle $B$ attached to $F$ and an oriented core $C$ of $B$.  Two $1$-handles with oriented cores $(B, C)$ and $(B', C')$ attached to $F$ are {\it equivalent} if they are ambient isotopic in $\R^4$ by an ambient isotopy of $\R^4$ keeping $F$ setwise fixed.  
(Here we assume that $C$ is mapped to $C'$ with respect to the orientations.)
\end{defn}

For a $1$-handle map $h$ attached to $F$, let $B$ be the image of $h$, which is a $1$-handle attached to $F$, and let $C$ be the core of $h$.  Using the core map $h^c: [0,1] \to \R^4$, we give an orientation to the core $C$.  Then we say that the $1$-handle with an oriented cord $(B, C)$ is {\it determined by $h$}.  

\begin{lem}\label{lem:twonotionsB}
For $1$-handles with oriented cores $(B, C)$ and $(B', C')$ attached to $F$, let $h$ and $h'$ be $1$-handle maps attached to $F$ determining $(B, C)$ and $(B', C')$, respectively.  $(B, C)$ and $(B', C')$ are equivalent if and only if $h$ and $h'$ are equivalent.  
\end{lem}

{\it Proof.} 
The if part is obvious. We prove the only if part.  It is sufficient to prove this in the case that $(B, C)=(B', C')$.  
Let $\partial h$ and $\partial h'$ be the restrictions of $h$ and $h'$ to $\partial ([0,1] \times B^2)$, respectively.  
Then $\partial h(\partial ([0,1] \times B^2))= \partial h'(\partial ([0,1] \times B^2)) = \partial B$, the initial point of $C$ is 
$\partial h ((0,0))= \partial h'((0,0))$ and  the terminal point of $C$ is $\partial h ((1,0))= \partial h'((1,0))$.  By a standard argument, so-called Alexander's trick, we see that $\partial h$ is ambient isotopic to $\partial h'$ in $\partial B$ keeping $h((0,0))$ and $h((1,0))$ fixed and keeping $F \cap B$ setwise fixed.  This ambient isotopy is extended to an ambient isotopy of $\R^4$ keeping $F$ setwise fixed.  So we may assume that $\partial h = \partial h'$.  By Alexander's trick, we may change $h$  so that $h=h'$, by an ambient isotopy of $B$ keeping $F \cap B$ setwise fixed, which is extended by an ambient isotopy of $\R^4$ keeping $F$ setwise fixed.  Thus $h$ is equivalent to $h'$.  \qed 

{\it Proof of Lemma~\ref{lem:twonotionsA}.}
Let $C$ and $C'$ be the oriented cores of $B$ and $B'$ such that $(B, C)$ and $(B', C')$ are determined by $h$ and $h'$, respectively.  
The $1$-handle $B$ is equivalent to $B'$ if and only if $(B, C)$ is equivalent to $(B', C')$ or $(B', -C')$.  
By Lemma~\ref{lem:twonotionsB}, $(B, C)$ is equivalent to $(B', C')$ if and only if $h$ is equivalent to $h'$, and 
$(B, C)$ is equivalent to $(B', -C')$ if and only if $h$ is equivalent to the reverse of $h'$. \qed

For a $1$-handle  $B$ attached to $F$, we denote by 
${\rm h}^1(F; B)$ the surface-knot 
$$(F -  {\rm int}(F \cap \partial B)) \cup (\partial B  -  {\rm int}(F \cap \partial B)),$$
which we call the surface-knot obtained from $F$ by {\it surgery} along $B$.  The surgery is also called a $1$-handle surgery or a hyperboloidal transformation (\cite{HK}).  

If $B$ and $B'$ are equivalent $1$-handles attached to $F$, then 
${\rm h}^1(F; B)$ and ${\rm h}^1(F; B')$ are ambient isotopic in $\R^4$.  

\begin{defn}
Assume that $F$ is an orientable surface-knot.  A $1$-handle  $B$ attached to $F$ is {\it orientable} (or {\it non-orientable}, resp.)  if ${\rm h}^1(F; B)$ is orientable (or non-orientable, resp.).     
\end{defn}

\section{Cords and $1$-handles attached to a surface-knot} 

Let $F$ be a surface-knot.  

\begin{defn}
A simple arc $C$ in $\R^4$ is a  {\it cord} attached to $F$ if $C \cap F= \partial C \cap F$ and if this  intersection consists of two distinct points of $F$.  An {\it oriented cord} is a cord with an orientation as a $1$-manifold.  
Two cords $C$ and $C'$ attached to $F$ are {\it equivalent} if they are ambient isotopic in $\R^4$ by an ambient isotopy of $\R^4$ keeping $F$ setwise fixed.  
\end{defn}

Let $F$ be a surface-knot and let $C$ be a cord attached to $F$.  Let $N(C)$ be a regular neighborhood of $C$ in $\R^4$, and put $F \cap N(C) =: M = M_- \cup M_+$, where $M_-$ and $M_+$ are disjoint 2-cells on $F$. (When $C$ is oriented, we assume that  the orientation of $C$ is from $M_-$ toward $M_+$.)  

Let $B$ be a $1$-handle attached to $F$ with core $C$.  We assume that $B$ is contained in ${\rm int }N(C)$.  
We denote by ${\rm h}^1(M; B)$ the surface  
$$(M -  {\rm int}(M \cap \partial B)) \cup (\partial B  -  {\rm int}(M \cap \partial B)),$$
which we call the surface obtained from $M$ by {\it surgery} along $B$.  Then 
${\rm h}^1(M; B) = {\rm h}^1(F; B) \cap N(C)$.  

\begin{defn}\label{def:compatibleo}
In the above situation, let $o$ be an orientation of $M = F \cap N(C)$.  
We say that $B$ is {\it compatible with $o$} if we can give an orientation to 
the surface ${\rm h}^1(M; B)$ such that the restriction to  $M -  {\rm int}(M \cap \partial B))$ 
of the orientation  coincides with that of  $o$.  Otherwise, we say that $B$ is {\it incompatible with $o$}.  
\end{defn}

\begin{lem}\label{lem:core1handleA}
For any cord $C$ attached to $F$ and for any orientation $o$ of $M= F \cap N(C)$, there exists a $1$-handle $B$ attached to $F$ with core $C$ contained in ${\rm int }N(C) \subset \R^4$ which is compatible with the orientation $o$.  Moreover, such a $1$-handle is unique up to ambient isotopy $(f_s \, | \, s \in [0,1])$ of $N(C)$ keeping $\partial N(C) \cup C$ pointwise fixed and $M$ setwise fixed.  
\end{lem}

{\it Proof. } By an ambient isotopy of $\R^4$, we move $F$, $C$, $N(C)$, $M= M_- \cup M_+$ so that 
$C = \{ (0,0,0, t) \in \R^4 \, | \, t \in [-1,1]\}$, 
$N(C) = \{ (x,y,z,t) \, | \, x^2 + y^2+ z^2 \leq 2, \, t \in [-2,2] \}$, 
$M_- = \{(x,y,0,-1) \, | \, x^2 + y^2 \leq 2\}$, 
$M_+ = \{(x,y,0,1) \, | \, x^2 + y^2 \leq 2\}$, and that 
the orientation of $o$ restricted to $M_-$ is opposite to that restricted to $M_+$.  It is sufficient to prove the lemma in the case where $F$, $C$, $N(C)$, $M$ and $o$ are in this situation.  

Let $N'(C)=  \{ (x,y,z,t) \, | \, x^2 + y^2+ z^2 \leq 1, \, t \in [-1-\epsilon,1+\epsilon] \}$ for a small positive number $\epsilon$, which is a smaller tubular neighborhood of $C$.  
 Let $B= \{ (x,y,0,t) \, | \, x^2 + y^2 \leq 1, \, t \in [-1,1] \}$.  It is a $1$-handle with core $C$ which is compatible with $o$.  Let $B'$ be another $1$-handle in ${\rm int }N(C)$ with core $C$ which is compatible with $o$. By an ambient isotopy $(f_s \, | \, s \in [0,1])$ of $N(C)$ keeping $\partial N(C) \cup C$ pointwise fixed and $M$ setwise fixed, we may assume that $B' =   \cup \{ X_t \times \{ t \}  \, | \, t \in [0,1] \}$ where $X_{-1}$ and $X_{1}$ are the standard 2-ball $B^2 \subset \R^2 \subset \R^3$, and for each $ t \in (-1, 1)$, $X_t$ is a unit 2-disk in $\R^3$ with center $0$.  We give $X_{-1}$ an orientation such that the orientation is the same with $o$ restricted to $M_- \cap B' = X_{-1} \times \{-1\}$.  For each $t \in (-1, 1]$, we give $X_t$ an orientation induced from the orientation of $X_{-1}$ continuously.  (Note that the orientation of $X_{1}$ is opposite to the orientation $o$ restricted to $M_- \cap B'  = X_{-1} \times \{1\}$, since $B'$ is compatible with $o$.) 
The one-parameter family $(X_t \, | \, t \in [-1,1])$ of oriented disks 
determines a family 
$(H_t \, | \, t \in [-1,1])$ of oriented 2-planes $H_t$  in $\R^3$ with $H_{-1} = H_{1} = \R^2$.  It induces a map $\theta: [-1, 1] \to G_{3,2}; t \mapsto H_t$ to the Grassmann manifold $G_{3,2}$ with $\theta(-1) = \theta(1)= \R^2$.  Since $G_{3,2}$ is homeomorphic to $S^2$, the loop $\theta$ is homotopic to the trivial map.  Hence by rotating the 3-balls $B^3 \times \{t\}$  for $t \in [-1,1]$ relative to $B^3 \times \{-1\} \cup B^3 \times \{1\}$, we can move $B'$ to $B$.  Using a collar neighborhood of $\partial N'(C)$ in $N(C)$, we may extend the rotations to an ambient isotopy $(g_s \, | \, s \in [0,1])$ of $N(C)$ keeping $\partial N(C) \cup C$ pointwise fixed and $M$ setwise fixed.  \qed 

\begin{lem}\label{lem:core1handleB}
For any cord $C$ attached to $F$ and for any orientation $o$ of $M= F \cap N(C)$, there exists a $1$-handle $B$ attached to $F$ with core $C$ contained in ${\rm int }N(C) \subset \R^4$ which is ``incompatible" with the orientation $o$.  Moreover, such a $1$-handle is unique up to ambient isotopy $(f_s \, | \, s \in [0,1])$ of $N(C)$ keeping $\partial N(C) \cup C$ pointwise fixed and $M$ setwise fixed.  
\end{lem}

{\it Proof. } By reversing the orientation of $M_+$ in the proof of Lemma~\ref{lem:core1handleA}, we see the result. \qed

We say that two cords $C$ and $C'$ attached to $F$ are {\it homotopic} if there is a homotopy $(C_s \, | \, s \in [0,1])$ consisting of arcs in $\R^4$ (possibly with self-intersection) with 
$C_0=C$ and $C_1=C'$ such that for each $s\in [0,1]$, $C_s \cap F = \partial C_s \cap F$ and the intersection consists of two distinct points of $F$.

\begin{lem}\label{lem:homotopyisotopyA}
Two cords $C$ and $C'$, with the same endpoints, $\partial C= \partial C'$, attached to $F$ are equivalent rel $\partial C$  if and only if they are homotopic rel $\partial C$.  
Two cords $C$ and $C'$ attached to $F$ are equivalent if and only if they are homotopic.  
\end{lem}

{\it Proof.} The former assertion follows from Theorem~4 of \cite{Hud}. The latter assertion is easily seen from the former. \qed

Let $C$ and $C'$ be cords attached to $F$ which are homotopic as cords attached to $F$ by a homotopy $(C_s \, | \, s \in [0,1])$ with $C_0=C$ and $C_1=C'$.  Let $o$ and $o'$ be  orientations of $M= F \cap N(C)$ and $M'= F \cap N(C')$, respectively.  For $s \in [0,1]$, let $o_s$ be the orientation of $M_s= F \cap N(C_s)$ induced from $o=o_0$ by the $1$-parameter family $(M_t \, | \, t \in [0, s])$.  If $o' = o_1$, then we say that the homotopy $(C_s \, | \, s \in [0,1])$  {\it respects  the orientations $o$ and $o'$}, and that $C$ and $C'$ are {\it homotopic with respect to the orientations $o$ and $o'$}.

\begin{defn}
Let $C$ and $C'$ be cords attached to $F$ which are equivalent as cords attached to $F$, and let $(f_s \, | \, s \in [0,1])$ be an ambient isotopy of $\R^4$ carrying $C$ to $C'$.  Consider the induced  homotopy $(C_s \, | \, s \in [0,1])$ with $C_s = f_s(C)$.  We say that $C$ and $C'$ are {\it equivalent with respect to the orientations $o$ and $o'$} if the homotopy $(C_s \, | \, s \in [0,1])$ respects $o$ and $o'$.  
\end{defn}

\begin{defn}
Let $C$ be an (oriented) cord attached to $F$ and $o$ be an orientation of $M= F \cap N(C)$.  
We call $o$ a {\it local orientation of $F$ at the endpoints of $C$}, and the pair 
$(C, o)$ an {\it oriented cord attached to $F$ with a local orientation at the endpoints}.  
Let $(C', o')$ be an oriented cord attached to $F$ with local orientation at the endpoints.  
We say that $(C, o)$ is {\it equivalent to $(C', o')$} if 
$C$ is equivalent to $C'$ with respect to $o$ and $o'$.  
\end{defn}

By the same reason of Lemma~\ref{lem:homotopyisotopyA}, we have the following lemma. 

\begin{lem}\label{lem:homotopyisotopyB}
Two cords $C$ and $C'$ attached to $F$ are equivalent with respect to $o$ and $o'$ (in other words, $(C, o)$ is equivalent to $(C', o')$) 
 if and only if they are homotopic with respect to $o$ and $o'$.  
\end{lem}

The following is a key lemma for classification of $1$-handles.  

\begin{lem}\label{lem:key}
Let $C$ and $C'$ be cords attached to $F$, and let $o$ and $o'$ be orientations of $M= F \cap N(C)$ and $M' = F \cap N(C')$.    Let $B$ and $B'$ be $1$-handles attached to $F$ with core $C$ and $C'$ such that they are compatible with $o$ and $o'$, respectively.  If $C$ and $C'$ are homotopic with respect to $o$ and $o'$, then $B$ and $B'$ are equivalent.  
\end{lem}

{\it Proof.} 
If $C$ and $C'$ are homotopic with respect to $o$ and $o'$, then by Lemma~\ref{lem:homotopyisotopyB}  they are equivalent with respect to $o$ and $o'$.  Let $(f_s \, | \, s \in [0,1])$ be an ambient isotopy of $\R^4$ carrying $C$ to $C'$ and keeping $F$ setwise fixed  such that the induced homotopy $(C_s \, | \, s \in [0,1])$ respects $o$ and $o'$, where $C_s:= f_s(C)$.  Let $f_1(B) = B''$.  By definition, $B$ is equivalent to $B''$.  Since $B''$ has core $C'$ and it is compatible with $o'$, by Lemma~\ref{lem:core1handleA}, we see that $B''$ is equivalent to $B'$.  Hence $B$ and $B'$ are equivalent.  \qed 

As a corollary, we have the following.   

\begin{cor}\label{cor:orientableA}
Let $F$ be an oriented surface-knot.  Let $B$ and $B'$ be orientable $1$-handles attached to $F$, and let $C$ and $C'$ be their cores.  
\begin{itemize}
\item[(1)] 
$B$ and $B'$ are equivalent if and only if $C$ and $C'$ are equivalent. 
\item[(2)] 
Suppose that $C$ and $C'$ are oriented. $(B, C)$ and $(B', C')$ are equivalent as $1$-handles with oriented cores if and only if $C$ and $C'$ are equivalent as oriented cords.  
\end{itemize}
\end{cor}

\begin{cor}\label{cor:orientableB}
Let $F$ be an oriented surface-knot.  Let $B$ and $B'$ be non-orientable $1$-handles attached to $F$, and let $C$ and $C'$ be their cores.  
\begin{itemize}
\item[(1)] 
$B$ and $B'$ are equivalent if and only if $C$ and $C'$ are equivalent. 
\item[(2)] 
Suppose that $C$ and $C'$ are oriented. $(B, C)$ and $(B', C')$ are equivalent as $1$-handles with oriented cores if and only if $C$ and $C'$ are equivalent as oriented cords.  
\end{itemize}
\end{cor}

\section{Classification of oriented cords} 

Let $F$ be a surface-knot, 
$N(F)$ be a tubular neighborhood of $F$  in  $\R^4$, and $E(F)$ be the exterior $\R^4 - {\rm int}N(F)$. 
The tubular neighborhood $N(F)$ is a $B^2$-bundle over $F$.  Let 
$p: N(F) \to F$ be the projection map.  A fiber $p^{-1}(y)$ $(y \in F)$ is called a {\it meridian disk} over $y$.  

Take a point $x$ in $\partial N(F) = \partial E(F)$ and put $G(F) := \pi_1(E(F), x)$, which is the {\it knot group} of $F$. 

\begin{sloppypar}
Let $\pi_1^+(\partial N(F), x)$ be the subgroup of $\pi_1(\partial N(F), x)$ consisting of all elements represented by loops in $\partial N(F)$ with base point $x$ such that their images under $p$ are orientation-preserving loops in $F$.  If $F$ is orientable, then $\pi_1^+(\partial N(F), x) = \pi_1(\partial N(F), x)$.  If $F$ is non-orientable, then 
$\pi_1^+(\partial N(F), x)$ is a subgroup of $\pi_1(\partial N(F), x)$ of index $2$.  
\end{sloppypar}

Let $P$ and $P^+$ denote subgroups of $G(F)$ that are the images of $\pi_1(\partial N(F), x)$ and 
$\pi_1^+(\partial N(F), x)$, respectively, under the inclusion-induced homomorphism  $i_\ast: \pi_1(\partial N(F), x) \to G(F)=\pi_1(E(F), x)$.  The subgroup $P$ is called the {\it peripheral subgroup} of $G(F)$, and we call 
$P^+$ the {\it positive peripheral subgroup}.

Let $C$ be an oriented cord attached to $F$.  We assume that $C \cap N(F)$ consists of two arcs each of which is  contained in a meridian disk.  The restriction of $C$ to $E(F)$ is an oriented simple arc in $E(F)$, which we denote by  $\overline{C}$.  Take a path $\alpha: [0,1] \to \partial N(F)$ such that $\alpha(0) = x $ and $\alpha(1)$ is the initial point of $\overline{C}$, and take a path $\beta : [0,1] \to \partial N(F)$ such that $\beta(0)= x$ and $\beta(1)$ is the terminal point of $\overline{C}$.  The composition $\alpha \overline{C} \beta^{-1}$ is a path in $E(F)$ with base point $x$, where we regard $\overline{C}$ as a path.   We have an element $[\alpha \overline{C} \beta^{-1}]$ of $G(F)$.  

\begin{defn}\label{defn:element}
In the above situation, we call the element $[\alpha \overline{C} \beta^{-1}]$ of $G(F)$ the {\it element determined from $C$ with $(\alpha, \beta)$}. 
\end{defn}

\begin{lem}[cf. \cite{Boy}]\label{lem:cordA}
The double coset $P[\alpha \overline{C} \beta^{-1}]P \in P \setminus   G(F) / P$ does not depend on a choice of $(\alpha, \beta)$.  
\end{lem}

{\it Proof.} 
Let $(\alpha', \beta')$ be another choice of paths.  
Then \begin{sloppypar} 
$P[\alpha \overline{C} \beta^{-1}]P$ $ = P[\alpha \alpha'^{-1} \alpha' \overline{C} \beta'^{-1} \beta' \beta^{-1}]P = 
P[\alpha' \overline{C} \beta'^{-1}]P$.  \qed 
\end{sloppypar}

For an oriented cord $C$ attached to $F$, 
we denote the double coset  $P[\alpha \overline{C} \beta^{-1}]P \in P \setminus  G(F) / P$ 
by $P(C)P$.  

The idea of the following theorem is essentially due to Boyle \cite{Boy}.  

\begin{thm}[cf. \cite{Boy}]\label{thm:cordA} 
Let $C$ and $C'$ be oriented cords attached to $F$.  
The cord $C$ is equivalent to $C'$ if and only if  $P(C)P = P(C')P$.  
\end{thm}

{\it Proof.} 
First we show that if $C$ is equivalent to $C'$, then $P(C)P = P(C')P$.  
Without loss of generality, we may assume that $C$ is ambient isotopic to $C'$ by 
an ambient isotopy of $\R^4$ keeping $F$ and $N(F)$ setwise fixed and keeping the base point $x$ fixed.  Let $(\alpha, \beta)$ be a pair of paths for $C$ as in Definition~\ref{defn:element}.  By the ambient isotopy of $\R^4$, let $(\alpha, \beta)$ be mapped to $(\alpha', \beta')$ and $\overline{C}$ be mapped to $\overline{C'}$.  Then 
$[\alpha \overline{C} \beta^{-1}] = [\alpha' \overline{C'} \beta'^{-1}]$ in $G(F)$.  Thus $P(C)P = P(C')P$.  

Suppose that $P(C)P = P(C')P$.  
Let $U$ be a regular neighborhood of $x$ in $\partial N(F)$, and let $\alpha_0$ and $\beta_0$ be short paths in $U$ with $\alpha_0(0)= \beta_0(0) = x$ and $\alpha_0(1) \neq \beta_0(1)$.  
Let $C$ and $C'$ be oriented cords attached to $F$ with $P(C)P = P(C')P$.  By moving $C$ and $C'$ up to equivalence, without loss of generality, we may assume that the starting points of $\overline{C}$ and $\overline{C'}$ are $\alpha_0(1)$ and the terminal points of $\overline{C}$ and $\overline{C'}$ are $\beta_0(1)$.  
Then 
$P[\alpha_0 \overline{C} \beta_0^{-1}]P = P[\alpha_0 \overline{C'} \beta_0^{-1}]P$, and hence $[\alpha_0 \overline{C} \beta_0^{-1}] = g[\alpha_0 \overline{C'} \beta_0^{-1}]g'$ in $G(F)$ for some elements $g, g' \in G(F)$.  
This implies that $\overline{C}$ is homotopic to $\overline{C'}$ in $E(F)$ after sliding the endpoints suitably, and we see that $C$ is homotopic to $C'$.  
By Lemma~\ref{lem:homotopyisotopyA}, $C$ is equivalent to $C'$. \qed 

\begin{thm}\label{thm:cordAA}
Let $\varphi$ be the map from the set of equivalence classes of oriented cords attached to $F$ to the double cosets $P \setminus   G(F) / P$ that sends the equivalence class of $C$ to $P(C)P$.  The map $\varphi$ is a bijection. 
\end{thm}

{\it Proof.} 
By Theorem~\ref{thm:cordA}, $\varphi$ is well defined and injective.  We show that $\varphi$ is surjective. 
Let $U$ be a regular neighborhood of $x$ in $\partial N(F)$, and let $\alpha_0$ and $\beta_0$ be short paths in $U$ with $\alpha_0(0)= \beta_0(0) = x$ and $\alpha_0(1) \neq \beta_0(1)$.  
Let $g$ be an element of $G(F)$.  There is a simple path $\gamma: [0,1] \to E(F)$ such that $g = [\alpha_0 \gamma \beta_0^{-1}]$ in $G(F)$.  Let $C$ be an oriented core attached to $F$ such that $\overline{C}$ is the image of $\gamma$.  
Then $P(C)P = PgP$.  Thus the map $\varphi$ is surjective.  \qed

Now we consider oriented cords attached to $F$ with local orientations of $F$ at the endpoints. 

Let $C$ be an oriented cord attached to $F$, and let $y_-$ and $y_+$ be the initial point and the terminal point of $C$, respectively.  
Let $M = M_- \cup M_+ = F \cap N(C)$, where $M_-$ is a 2-cell in $F$ containing $y_-$ and $M_+$ is a 2-cell containing $y_+$.  We denote by $y$ the image $p(x)$ of $x$, and let $U_y$ be a regular neighborhood of $y$ in $F$.  
Let $o_y$ be an orientation of $U_y$, and   
let $o$ be an orientation of $M$.  

Let $(\alpha, \beta)$ be a pair of paths for $C$ as in Definition~\ref{defn:element}.  Then $(p \alpha, p \beta)$ is a pair of paths in $F$ with $(p \alpha)(0)=(p \beta)(0) =y$, $(p \alpha)(1) = y_-$ and $(p \beta)(1)= y_+$.  

\begin{defn}\label{def:compatibleoyo}
We say that $(\alpha, \beta)$ is {\it compatible with $o_y$ and $o$} if 
the local orientation of $F$ at $y_-$ determined from $o$ coincides with the local orientation at $y_-$ obtained from the local orientation $o_y$ at $y$ by translating along $p \alpha$ and if   
the local orientation of $F$ at $y_+$ determined from $o$ coincides with the local orientation at $y_+$ obtained from the local orientation $o_y$ at $y$ by translating along $p \beta$.  Otherwise, we say that it is {\it incompatible with $o_y$ and  $o$}.  
\end{defn}

\begin{lem}\label{lem:cordB}
Let $o_y$ be an orientation of $U_y$.  
Let $C$ be an oriented cord attached to $F$, and $o$ be an orientation of $M= F \cap N(C)$.  
Let $(\alpha, \beta)$ be a pair of paths as in Definition~\ref{defn:element}.  
Suppose that $(\alpha, \beta)$ is compatible with $o_y$ and $o$.  
The double coset $P^+[\alpha \overline{C} \beta^{-1}]P^+ \in P^+ \setminus   G(F) / P^+$ does not depend on a choice of $(\alpha, \beta)$.  
\end{lem}

{\it Proof.} 
Let $(\alpha', \beta')$ be another choice of paths that is compatible with $o_y$ and $o$.  
Since $(p \alpha)(p \alpha'^{-1})$ is an orientation-preserving loop in $F$, we have that $[\alpha \alpha'^{-1}] \in G(F)$ belongs to $P^+$.  Similarly, $[\beta' \beta^{-1}]$ belongs to $P^+$.  Then 
\begin{sloppypar}
$P^+[\alpha \overline{C} \beta^{-1}]P^+$ $ = P^+[\alpha \alpha'^{-1} \alpha' \overline{C} \beta'^{-1} \beta' \beta^{-1}]P^+ = 
P^+[\alpha' \overline{C} \beta'^{-1}]P^+$.  \qed 
\end{sloppypar}

In the situation of Lemma~\ref{lem:cordB}, we denote the double coset $P^+[\alpha \overline{C} \beta^{-1}]P^+ 
\in P^+ G(F) P^+$ by $P^+ (C, o_y, o) P^+$.

\begin{thm}\label{thm:cordB} 
Let $o_y$ be an orientation of $U_y$.  
Let $(C, o)$ and $(C', o')$ oriented cords attached to $F$ with local orientations at the endpoints.  
$(C, o)$ is equivalent to $(C', o')$ if and only if  
$P^+(C, o_y, o)P^+ = P^+(C', o_y, o')P^+$.   
\end{thm}

The proof is analogous to the proof of Theorem~\ref{thm:cordA}, which is left to the reader. 

\begin{thm}\label{thm:cordBB} 
Let $o_y$ be an orientation of $U_y$.  
Let $\psi$ be the map from the set of equavalence classes of oriented cords attached to $F$ with local orientations at the endpoints 
to the double cosets $P^+ \setminus   G(F) / P^+$ that sends the equivalence class of $(C, o)$  to $P^+(C, o_y, o)P^+$.  The map $\psi$ is a bijection.  
\end{thm}

{\it Proof.} (The proof is analogous to the proof of Theorem~\ref{thm:cordAA}.) 
By Theorem~\ref{thm:cordB}, $\psi$ is well defined and injective.  We show that $\psi$ is surjective. 
Let $U$ be a regular neighborhood of $x$ in $\partial N(F)$, and let $\alpha_0$ and $\beta_0$ be short paths in $U$ with $\alpha_0(0)= \beta_0(0) = x$ and $\alpha_0(1) \neq \beta_0(1)$.  
Let $g$ be an element of $G(F)$.  There is a simple path $\gamma: [0,1] \to E(F)$ such that $g = [\alpha_0 \gamma \beta_0^{-1}]$ in $G(F)$.  Let $C$ be an oriented core attached to $F$ such that $\overline{C}$ is the image of $\gamma$.  Let $o$ be a local orientation of $F$ at the endpoints of $C$ such that it is obtained from $o_y$ by translating $o_y$ along $p \alpha_0$ and $p \beta_0$.  
Then $\psi(C, o) = P^+(C, o_y, o)P^+ = P^+ g P^+$.  Thus the map $\psi$ is surjective.  \qed

\section{Classification of  $1$-handles in Case 1} 

In this section we consider Case 1:  $F$ is oriented and $1$-handles are orientable.  

First we classify $1$-handles with oriented cores.  This case is due to Boyle \cite{Boy}.  

Let $F$ be an oriented surface-knot, $G(F)$ be the knot group $\pi_1(E(F), x)$, and $P$ be the peripheral subgroup.  

Let $(B,C)$ be an orientable $1$-handle with an oriented core attached to $F$.  
Let $P(C)P$ be the double coset $[\alpha \overline{C} \beta^{-1}] \in P \setminus G(F) /P
$ where  $\overline{C} = C \cap E(F)$ and 
$(\alpha, \beta)$ is a pair of paths for $C$ as in Definition~\ref{defn:element}.  

\begin{thm}[Boyle \cite{Boy}]\label{thm:case1A}
Two orientable $1$-handles with oriented cores $(B,C)$ and $(B',C')$ attached to $F$ are equivalent if and only if $P(C)P=P(C')P$.  
Moreover, a map sending the equivalence class of $(B,C)$ to $P(C)P$ is a bijection from the set of equivalence classes of orientable $1$-handles with oriented cores attached to $F$ to the double cosets $P \setminus   G(F) / P$.  
\end{thm}

{\it Proof.} 
By Corollary~\ref{cor:orientableA}, the equivalence class of $(B, C)$ corresponds to the equivalence class of the oriented core $C$.  By Theorem~\ref{thm:cordAA}, we have the result. \qed 

\begin{defn}
Let $B$ be an orientable $1$-handle attached to $F$.  
Define $P(B)P $ by an unordered pair $\{ P(C)P, P(-C)P \}$, where $C$ is  an oriented core $C$ of $B$.  
\end{defn}

Note that $P(B)P $ does not depend on a choice of an oriented core of $B$.  

\begin{thm}\label{thm:case1B}
Two orientable $1$-handles $B$ and $B'$ attached to $F$ are equivalent if and only if $P(B)P= P(B')P$.  
\end{thm}

{\it Proof.} Let $C$ and $C'$ be oriented cores of $B$ and $B'$, respectively.  
Then $B$ and $B'$ are equivalent if and only if $(B, C)$ is equivalent to $(B', \epsilon C')$ for some $\epsilon \in \{ +1, -1\}$. 
The latter statement holds if and only if $(B, C)$ is equivalent to $(B', \epsilon C')$ and $(B, -C)$ is equivalent to $(B', -\epsilon C')$, and hence by Theorem~\ref{thm:case1A},  if and only if $P(C)P= P(\epsilon C')P$ and $P(-C)P= P(-\epsilon C')P$. The last statement holds if and only if $P(B)P = P(B')P$. \qed 

For a set $X$ and a positive integer $n$, the {\it unordered $n$-fold product of $X$} means the set of unordered $n$-tuples of elements $\{ x_1, \dots, x_n \}$ of $X$.

Note that $P(B)P $ is an element of the unordered $2$-fold product of 
$P \setminus   G(F) / P$.  

By Theorem~\ref{thm:case1B}, 
a map sending the equivalence class of $B$ to $P(B)P$  from the set of equivalence classes of orientable $1$-handles attached to $F$ to the unordered $2$-fold product of $P \setminus   G(F) / P$ is well defined and injective.  
However this map is not surjective in general.  The image of this map is characterized as follows.

\begin{prop}\label{prop:case1image}
The image of the map sending the equivalence class of $B$ to $P(B)P$  from the set of equivalence classes of orientable $1$-handles attached to $F$ to the unordered $2$-fold product of $P \setminus   G(F) / P$ is the subset consisting of the elements $\{ PgP, Pg^{-1}P \}$ for all $g \in G(F)$.  
\end{prop}

{\it Proof.} Let $B$ be an orientable $1$-handle attached to $F$, and $C$ an oriented core.  
By definition, $P(C)P = P[\alpha \overline{C} \beta^{-1}]P$ as before.  Put $g = [\alpha \overline{C} \beta^{-1}] \in G(F)$.  
Then $P(-C)P = P[\beta \overline{C}^{-1} \alpha^{-1}]P = Pg^{-1}P$.  Hence $P(B)P = \{ PgP, Pg^{-1}P \}$ for some $g \in G(F)$.  Conversely for any $g \in G(F)$, there is an orientable $1$-handle with an oriented core $C$ such that $P(C)P = PgP$.  Then $P(B)P = \{ PgP, Pg^{-1}P \}$. \qed 

\begin{cor}\label{for:case1notsurjective}
Let $F$ be an oriented surface-knot with $P \neq G(F)$.  Then the map sending the equivalence class of $B$ to $P(B)P$  from the set of equivalence classes of orientable $1$-handles attached to $F$ to the unordered $2$-fold product of $P \setminus   G(F) / P$ is not surjective.  
\end{cor}

{\it Proof.} 
Let $g$ be an element of  $G(F) - P$.  
We show that $\{ PgP, P1P \}$ is not obtained from any $1$-handle.  
Assume that $\{ PgP, P1P \}= \{ Pg'P, Pg'^{-1}P \}$ for some $g' \in G(F)$. 
Replacing $g'$ with $g'^{-1}$ if it is necessary, 
we may assume that $PgP = Pg'P$.  Since $g \notin P$,  we have $g' \notin P$.  On the other hand, $P1P = Pg'^{-1}P$ implies that $g' \in P$. This is a contradiction. \qed

For example, every non-trivial $2$-knot satisfies that $P \neq G(F)$.

\section{Classification of $1$-handles in Case 2} 

In this section we consider Case 2: $F$ is oriented and $1$-handles are non-orientable.  

This case is completely analogus to Case~1.  
Let $F$ be an oriented surface-knot, $G(F)$ be the knot group $\pi_1(E(F), x)$, and $P$ be the peripheral subgroup.  

Let $(B,C)$ be a non-orientable $1$-handle with an oriented core attached to $F$.  
Let $P(C)P$ be the double coset $[\alpha \overline{C} \beta^{-1}] \in P \setminus G(F) /P
$ where  $\overline{C} = C \cap E(F)$ and 
$(\alpha, \beta)$ is a pair of paths for $C$ as in Definition~\ref{defn:element}.  

\begin{thm}\label{thm:case2A}
Two non-orientable $1$-handles with oriented cores $(B,C)$ and $(B',C')$ attached to $F$ are equivalent if and only if $P(C)P=P(C')P$.  
Moreover, a map sending the equivalence class of $(B,C)$ to $P(C)P$ is a bijection from the set of equivalence classes of non-orientable $1$-handles with oriented cores attached to $F$ to the double cosets $P \setminus   G(F) / P$.  
\end{thm}

{\it Proof.} 
By Corollary~\ref{cor:orientableB}, the equivalence class of $(B, C)$ corresponds to the equivalence class of the oriented core $C$.  By Theorem~\ref{thm:cordAA}, we have the result. \qed 

\begin{defn}
Let $B$ be a non-orientable $1$-handle attached to $F$.  
Define $P(B)P $ by an unordered pair $\{ P(C)P, P(-C)P \}$, where $C$ is  an oriented core $C$ of $B$.  
\end{defn}

Note that $P(B)P $ does not depend on a choice of an oriented core of $B$.  

\begin{thm}\label{thm:case2B}
Two non-orientable $1$-handles $B$ and $B'$ attached to $F$ are equivalent if and only if $P(B)P= P(B')P$.  
\end{thm}

{\it Proof.} The proof is the same with the proof of Theorem~\ref{thm:case1B}, where we use Theorem~\ref{thm:case2A} instead of Theorem~\ref{thm:case1A}.  \qed 

By Theorem~\ref{thm:case2B}, 
a map sending the equivalence class of $B$ to $P(B)P$  from the set of equivalence classes of non-orientable $1$-handles attached to $F$ to the unordered $2$-fold product of $P \setminus   G(F) / P$ is well defined and injective.  
This map is not surjective in general.  The image of this map is exactly the same with the subset given in Proposition~\ref{prop:case1image}.  An analogous statement to Corollary~\ref{for:case1notsurjective} is also valid for non-orientable $1$-handles.

\section{Classification of  $1$-handles in Case 3} 

In this section we consider Case 3: $F$ is non-orientable.  

Let $F$ be a non-oriented surface-knot, $G(F)$ be the knot group $\pi_1(E(F), x)$, and $P^+$ be the positive peripheral subgroup.  Let $o_y$ be an orientation of a regular neighborhood $U_y$ of $y= p(x)$.  

First we classify $1$-handles with oriented cores. 

Let $(B,C)$ be a $1$-handle with an oriented core attached to $F$.  
Let $o$ be an orientation of $M= F \cap N(C)$ such that the $1$-handle   $B$ is compatible with $o$ (Definition~\ref{def:compatibleo}). 
Let $(\alpha, \beta)$ be a pair of paths for $C$ as in Definition~\ref{defn:element} such that it is compatible with $o_y$ and $o$ (Definition~\ref{def:compatibleoyo}).  
Let $P^+(C, o_y, o)P^+$ be the double coset $[\alpha \overline{C} \beta^{-1}] \in P^+ \setminus G(F) /P^+
$ where  $\overline{C} = C \cap E(F)$.  

By $-o_y$ and $-o$, we denote the reversed orientations of $o_y$ and $o$, respectively.

\begin{lem}\label{lem:case3pre}
In the situation above, we have the following. 
\begin{itemize}
\item[(1)] 
The $1$-handle $B$ is compatible with an orientation $o'$ of $M$ if and only if  $o'=o$ or  $o'= -o$.  
\item[(2)] 
$P^+(C, -o_y, o)P^+ = P^+(C, o_y, -o)P^+$, and \\ 
$P^+(C, -o_y, -o)P^+ = P^+(C, o_y, o)P^+$.  
\end{itemize}
\end{lem}

{\it Proof.} 
(1) Since $M = F \cap N(C) = M_- \cup M_+$, there are four orientations of $M$.  Let $o | M_-$ and $o | M_+$ be the orientations of $M_-$ and $M_+$ that are restrictions of $o$.  Then the four orientations are 
$o = (o | M_-, o | M_+), -o= (-o | M_-, -o | M_+), (o | M_-, -o | M_+)$ and $(-o | M_-, o | M_+)$.  If $B$ is compatible with $o$ then it is compatible with $-o$ and it is incompatible with the other two.  

(2) It is obvious by the definition of $P^+(C, o_y, o)P^+$.  \qed 

\begin{defn}
Let $(B, C)$ be a $1$-handle with an oriented core attached to a non-orientable surface-knot $F$.  In the situation above, 
we define $P^+(B, C)P^+$ by an unordered pair $\{ P^+(C, o_y, o)P^+, P^+(C, o_y, -o)P^+ \}$, which is an element of the unordered $2$-fold product of $P^+ \setminus G(F) /P^+$.  
\end{defn}

By Lemma~\ref{lem:case3pre}, $P^+(B, C)P^+$ does not depend on a choice of $o_y$ and $o$.  

\begin{thm}\label{thm:case3A}
Let $F$ be a non-orientable surface-knot. 
Two $1$-handle with oriented cores  $(B,C)$ and $(B', C')$ attached to $F$ are equivalent if and only 
$P^+(B, C)P^+ = P^+(B', C')P^+$.  
\end{thm}

{\it Proof.} 
Let $o$ and $o'$ be orientations of $M = F \cap N(C)$ and $M'= F \cap N(C')$ such that $B$ and $B'$ are compatible with $o$ and $o'$, respectively.  

Note that $(B,C)$ and $(B', C')$ are equivalent if and only if $(C, o)$ is equivalent to $(C', \epsilon o')$ for some $\epsilon \in \{ +1, -1\}$, and $(C, -o)$ is equivalent to $(C', -\epsilon o')$.  
By Theorem~\ref{thm:cordB}, this condition is equivalent to that $P^+(C, o_y, o)P^+= P^+(C', o_y, \epsilon o')P^+$ and $P^+(C, o_y, -o)P^+= P^+(C', o_y, -\epsilon o')P^+$.  It is equivalent to that $P^+(B, C)P^+ = P^+(B', C')P^+$. \qed 

\begin{sloppypar}

\begin{defn}\label{case3mapA}
Let $\theta$ be a 
map sending the equivalence class of $(B, C)$ to $P^+(B, C)P^+$ from the set of equivalence classes of $1$-handles with oriented cores attached to $F$ to the unordered $2$-fold product of  $P^+ \setminus   G(F) / P^+$.  
\end{defn}

\end{sloppypar}

By Theorem~\ref{thm:case3A}, the map $\theta$ is well defined and injective.  
In general, it is not surjective.  

\begin{lem}\label{lem:case3imageA1}
Let $m$ be an element of $\pi_1(\partial N(F), x) - \pi_1^+(\partial N(F), x)$.  
The image of the map $\theta$ is the subset consisting of the elements $\{ P^+ g P^+, P^+ i_\ast(m) \, g \, i_\ast(m) P^+ \}$ for all $g \in G(F)$.  
\end{lem}

{\it Proof.} We may fix an orientation $o_y$.   Let $\nu$ be a loop in $\partial N(F)$ with base point $x$ representing $m$.  
Let $(B, C) $ be a $1$-handle with an oriented core attached to $F$.   Let $o$ be an orientation of $M = F \cap N(C)$ such that  $B$ is compatible with $o$.  (Then $B$ is also compatible with $-o$.) 
Let $(\alpha, \beta)$ be a pair of paths in $N(C)$ for $C$ as in Definition~\ref{defn:element} such that it is compatible with $o_y$ and $o$.  
Then by definition, $P^+(C, o_y, o)P^+ = P^+ [\alpha \overline{C} \beta^{-1}] P^+$, 
where  $\overline{C} = C \cap E(F)$ and we regard $\overline{C}$ as a simple path.  
Note that $(\nu \alpha, \nu^{-1} \beta)$ is a pair of paths in $N(F)$ for $C$  such that it is compatible with $o_y$ and $-o$.  
Thus $P^+(C, o_y, -o)P^+ = P^+ [\nu \alpha \overline{C} \beta^{-1} \nu]P^+$.  
When we put $g = [\alpha \overline{C} \beta^{-1}] \in G(F)$, we have 
$P^+(C, o_y, o)P^+ = P^+ g P^+$ and $P^+(C, o_y, -o)P^+ = P^+ i_\ast(m) \, g \, i_\ast(m) P^+$.  
Hence $P^+(B, C)P^+ =  \{ P^+ g P^+, P^+ i_\ast(m) \, g \, i_\ast(m) P^+ \}$.  

Conversely for any $g \in G(F)$, there is a $1$-handle with an oriented core $(B,C)$ attached to $F$  such that 
$P^+(C, o_y, o)P^+ = P^+ g P^+$.  This is verified by the same argument in the proof of Theorem~\ref{thm:cordBB}.  
Then, as shown above, we see that $P^+(B, C)P^+ =  \{ P^+ g P^+, P^+ i_\ast(m) \, g \, i_\ast(m) P^+ \}$.  \qed 

\begin{prop}\label{prop:case3imageA2}
The image of the map $\theta$ is characterized as follows:  
\begin{itemize}
\item[(1)] If $P^+ \neq P$, then 
the image  consists of the elements 
$$\{ P^+ g P^+, P^+ n g n P^+ \}$$ for all $g \in G(F)$, 
where $n$ is an element of $P - P^+$. 
\item[(2)] 
If $P^+ =  P$, then 
the image  consists of the elements 
$$\{ P^+ g P^+, P^+ g P^+ \}$$ for all $g \in G(F)$. 
\end{itemize}
\end{prop}

{\it Proof.} 
(1) Suppose that $P^+ \neq P$ and let $n \in P - P^+$.  Take 
an element $m \in \pi_1(\partial N(F), x) -  \pi_1^+(\partial N(F), x)$ with $i_\ast (m)= n$.  
By Lemma~\ref{lem:case3imageA1}, we have the result. 

(2) Suppose that $P^+ = P$.  There is an element $m \in \pi_1(\partial N(F), x) -  \pi_1^+(\partial N(F), x)$ with $i_\ast (m) =1$. By Lemma~\ref{lem:case3imageA1}, we have the result.  \qed

\begin{cor}\label{for:case3notsurjective1}
Let $F$ be a non-orientable surface-knot with $P^+ \neq P$.  The map $\theta$ is not surjective.  
\end{cor}

{\it Proof.} 
Let $n$ be an element of $P -  P^+$.  
We show that  $\{ P^+ nP^+, P^+ 1P^+  \}$ is not in the image of $\theta$.  
Assume that $\{ P^+ nP^+, P^+ 1P^+  \}$ is in the image of $\theta$.  
Then by Proposition~\ref{prop:case3imageA2} (1), 
$\{ P^+ nP^+, P^+ 1P^+ \}= \{ P^+ g P^+, P^+ n g n P^+ \}$ for some $g \in G(F)$. 

(i) Suppose that $P^+ nP^+ = P^+ g P^+$ and $P^+ 1P^+ = P^+ n g n P^+$.  
Since $n \in P -  P^+$, $P^+ nP^+ = P^+ g P^+$ implies that $ g \in P - P^+$.  On the other hand,  $P^+ 1P^+ = P^+ n g n P^+$ implies that $g \in P^+$.  This is a contradicition.  

(ii) Suppose that $P^+ nP^+ = P^+ n g n P^+$ and $P^+ 1P^+ = P^+ g P^+$.  
Since $n \in P -  P^+$, $P^+ nP^+ = P^+ n g n P^+$ implies that $ g \in P - P^+$.  On the other hand, $P^+ 1P^+ = P^+ g P^+$ implies that $g \in P^+$.  This is a contradicition.  \qed

A surface-knot $F$ is said to be {\it incompressible} if the inclusion-induced homomorphism 
$i_\ast: \pi_1(\partial N(F), x) \to G(F)$ is injective.  A method of constructing incompressible Klein bottles in $R^4$ is given in \cite{Kam90}.  
(The method in \cite{Kam90} relied on the existence of incompressible tori in $\R^4$, which is shown in \cite{Asa, Lit}.) 

Incompressible non-orientable surface-knots satisfy that $P^+ \neq P$.  

\begin{cor}\label{for:case3notsurjective2}
Let $F$ be a non-orientable surface-knot with $P^+ = P \neq G(F)$.  The map $\theta$ is not surjective.  
\end{cor}

{\it Proof.} 
Let $g$ be an element of $G(F) - P^+$.  
We show that  $\{ P^+ gP^+, P^+ 1P^+  \}$ is not in the image of $\theta$.  
Assume that $\{ P^+ gP^+, P^+ 1P^+  \}$ is in the image of $\theta$.  
Then by Proposition~\ref{prop:case3imageA2} (2), 
$\{ P^+ gP^+, P^+ 1P^+ \}= \{ P^+ g' P^+, P^+  g' P^+ \}$ for some $g' \in G(F)$. 
This implies that $g \in P^+$.  This contradicts to $g \in G(F) - P^+$. \qed 

If $F$ is a non-orientable surface-knot  obtained from a connected sum of a standard projective plane in $\R^4$ (cf. \cite{HK}) and a surface-knot, then $F$ satisfies that $P^+ = P$.  Using this, one can obtain a lot of examples of non-orientable surface-knots with $P^+ = P \neq G(F)$.

Now we consider $1$-handles attached to $F$.  

\begin{defn}
Let $B$ be a  $1$-handle attached to $F$.  
Define $P^+(B)P^+ $ by an unordered pair of unordered pairs $\{ P^+(B, C)P^+, P^+(B, -C)P^+ \}$, where $C$ is  an oriented core $C$ of $B$.  \end{defn}

Note that $P^+(B)P^+ $ does not depend on a choice of an oriented core of $B$.  

\begin{thm}\label{thm:case3B}
Two $1$-handles $B$ and $B'$ attached to $F$ are equivalent if and only if $P^+(B)P^+= P^+(B')P^+$.  
\end{thm}

{\it Proof.} Let $C$ and $C'$ be oriented cores of $B$ and $B'$, respectively.  
Then $B$ and $B'$ are equivalent if and only if $(B, C)$ is equivalent to $(B', \epsilon C')$ for some $\epsilon \in \{ +1, -1\}$. 
The latter statement holds if and only if $(B, C)$ is equivalent to $(B', \epsilon C')$ and $(B, -C)$ is equivalent to $(B', -\epsilon C')$, and hence by Theorem~\ref{thm:case3A},  if and only if $P^+(B,C)P^+= P^+(B', \epsilon C')P^+$ and $P^+(B, -C)P^+= P^+(B', -\epsilon C')P^+$. The last statement holds if and only if $P^+(B)P^+= P^+(B')P^+$. \qed 

Let $J$ be the image of the map $\theta$ defined in Definition~\ref{case3mapA}, which is characterized in Lemma~\ref{lem:case3imageA1} and Proposition~\ref{prop:case3imageA2}.  

\begin{defn}\label{case3mapB}
Let $\Theta$ be a map sending the equivalence class of $B$ to $P^+(B)P^+$  from the set of equivalence classes of $1$-handles attached to $F$ to the unordered $2$-fold product of $J$.  
\end{defn}

By Theorem~\ref{thm:case3B},  the map $\Theta$ is well defined and injective.   

By Lemma~\ref{lem:case3imageA1} and Proposition~\ref{prop:case3imageA2}, we have the following.  

\begin{lem}\label{lem:case3imageB1}
Let $m$ be an element of $\pi_1(\partial N(F), x) -  \pi_1^+(\partial N(F), x)$.   
The image of the map $\Theta$ is the subset consisting of the elements 
$$ \{    \{ P^+ g P^+, P^+ i_\ast (m) \, g \, i_\ast (m) P^+ \}, \{ P^+ g^{-1} P^+, P^+ i_\ast (m) \, g^{-1} \, i_\ast (m) P^+ \} \}$$
for all $g \in G(F)$.  
\end{lem}

{\it Proof.} 
This follows from Lemma~\ref{lem:case3imageA1} by a similar argument as in the proof of Proposition~\ref{prop:case1image}. \qed  

\begin{lem}\label{lem:case3imageB2}
The image of the map $\Theta$ is characterized as follows:  
\begin{itemize}
\item[(1)] If $P^+ \neq P$, then 
the image consists of the elements 
$$ \{    \{ P^+ g P^+, P^+ n g n P^+ \}, \{ P^+ g^{-1} P^+, P^+ n g^{-1} n P^+ \} \}$$ for all $g \in G(F)$, 
where $n$ is an element of $P - P^+$. 
\item[(2)] 
If $P^+ =  P$, then 
the image consists of the elements 
$$ \{    \{ P^+ g P^+, P^+ g P^+ \}, \{ P^+ g^{-1} P^+, P^+  g^{-1}  P^+ \} \}$$ for all $g \in G(F)$. 
\end{itemize}
\end{lem}

{\it Proof.} 
This follows from Proposition~\ref{prop:case3imageA2} by a similar argument as in the proof of Proposition~\ref{prop:case1image}. \qed

\vspace{0.4in}

\end{document}